\documentclass[twoside,12pt, leqno]{article}
\usepackage{amsmath,amscd,amsthm,amssymb,amsxtra,latexsym,epsfig,epic,graphics}
\usepackage[matrix,arrow,curve]{xy}
\usepackage{graphicx}
\usepackage{diagrams}
\def\antiddot{\mathinner{\mkern1mu\raise1pt\vbox{\kern7pt\hbox{.}}\mkern2mu
        \raise4pt\hbox{.}\mkern2mu\raise7pt\hbox{.}\mkern1mu}}

\newcommand{\KK}{{\mathbb K}}
\newcommand{\PP}{{\mathbb P}}
\newcommand{\QQ}{{\mathbb Q}}
\newcommand{\RR}{{\mathbb R}}

\newcommand{\ZZ}{{\mathbb Z}}

\newcommand{\Ext}{{\rm{Ext}}}

\newcommand{\s}{\mathcal}

\newcommand{\cO}{{\s O}}
\newcommand{\cF}{{\s F}}

\newcommand{\cI}{{\s I}}

\newcommand{\punkt}{\hspace{-.3ex}\raise.15ex\hbox to1ex{\Huge.}}

\def \fix#1 {{\hfill\break \bf (( #1 ))\hfill\break}}

\DeclareMathOperator{\Hom}{Hom}

\DeclareMathOperator{\Tor}{Tor}

\newtheorem{theorem}{Theorem}[section]
\newtheorem{lemma}[theorem]{Lemma}
\newtheorem{proposition}[theorem]{Proposition}
\newtheorem{corollary}[theorem]{Corollary}

\theoremstyle{definition}

\newtheorem{example}[theorem]{Example}

\newtheorem{algorithm}[theorem]{Algorithm}

\def\QQ{{\mathbb Q}}

\def\cF{{\cal F}}

\makeatletter
\def\Ddots{\mathinner{\mkern1mu\raise\p@
\vbox{\kern7\p@\hbox{.}}\mkern2mu
\raise4\p@\hbox{.}\mkern2mu\raise7\p@\hbox{.}\mkern1mu}}
\makeatother

\makeatletter
\def\Ddots{\mathinner{\mkern1mu\raise\p@
\vbox{\kern7\p@\hbox{.}}\mkern2mu
\raise4\p@\hbox{.}\mkern2mu\raise7\p@\hbox{.}\mkern1mu}}
\makeatother

\usepackage{times}
\newdimen\x \x=12pt

\usepackage{color}

\definecolor{Gray}{gray}{.75}

\def\d_#1{\underline{d_{#1}}}
\def\z_#1{\underline{z_{#1}}}
\def\r_#1{\underline{\lower.5pt\vbox{}\smash{r_{#1}}}}

\newdimen\eps\eps=4pt
\def\grayblock #1 #2 {\vrule width#1\x height#2\x depth0pt}
\def\graystripe #1 #2 {\raise #2\x\hbox{\vrule width#1\x height 0pt depth\x}}

\def\step #1 #2 #3 {\raise #2\x\hbox{\vrule width#1\x height 0.2pt depth0.2pt}%
  \hskip-.2pt\vrule height#2\x depth-#3\x\hskip-.2pt}

\def\backtrack{\hskip-27\x}

\def\plop $#1$ at (#2,#3){\raise #3\x
  \hbox to 0pt{\hskip #2\x\raise\eps\hbox to\x{$\hss#1\hss$}\hss}\ignorespaces}
\def\lowerplop $#1$ at (#2,#3){\raise #3\x
  \hbox to 0pt{\hskip #2\x\raise3pt\hbox to\x{$\hss#1\hss$}\hss}\ignorespaces}

\def\Plop $#1$ at (#2,#3){\raise #3\x
  \hbox to 0pt{\hskip #2\x\vbox{
\vskip-.2pt\hrule\vskip-.2pt\hbox{%
\hskip-.2pt\vrule height\x depth0pt\hskip-.2pt
\raise \eps\hbox to\x{$\hss#1\hss$}%
\hskip-.2pt\vrule height\x depth0pt\hskip-.2pt}%
\vskip-.2pt\hrule\vskip-.2pt}\hss}\ignorespaces}

\def\jscale{%
\hskip28\x \step 0 8 0 \hskip-28\x 
\lowerplop $0$ at (28,0)
\lowerplop $1$ at (28,1)
\lowerplop $2$ at (28,2)
\lowerplop $3$ at (28,3)
\lowerplop $4$ at (28,4)
\lowerplop $5$ at (28,5)
\lowerplop $6$ at (28,6)
\lowerplop $i$ at (28,7.5)
}

\def\jscaleo{%
\hskip26\x \step 0 10 0 \hskip-26\x 
\lowerplop $0$ at (26,0)
\lowerplop $1$ at (26,1)
\lowerplop $2$ at (26,2)
\lowerplop $3$ at (26,3)
\lowerplop $4$ at (26,4)
\lowerplop $5$ at (26,5)
\lowerplop $6$ at (26,6)
\lowerplop $7$ at (26,7)
\lowerplop $8$ at (26,8)
\lowerplop $i$ at (26,9.5)
}

\def\background{%
\grayblock 7 7  
\grayblock 5 4
\grayblock 7 3
\grayblock 3 2
\grayblock 5 1
\endgroup
\backtrack
\step 7 7 4     
\step 5 4 3
\step 7 3 2
\step 3 2 1
\step 5 1 1     
\backtrack
\step 27 0 0    
\backtrack
}
\def\backgroundo{%
\grayblock 27 8  
\endgroup
\backtrack
\step 7 8 8     
\step 7 8 8
\step 7 8 8
\step 3 8 8
\step 4 8 8     
\backtrack
\step 27 0 0    
\backtrack
}

\date{}
\title{Cohomology of  Coherent Sheaves and Series of Supernatural Bundles}
\author{David Eisenbud and Frank-Olaf Schreyer}

\begin{document}

\maketitle

\begin{abstract}

We show that
the cohomology
table of any coherent sheaf  on projective space is a
convergent---but possibly infinite---sum of positive real multiples of the cohomology tables of
what we call \emph{supernatural}
sheaves. 

\end{abstract}

\section*{Introduction}

Let $\KK$ be a field, and let $\cF$ be a coherent sheaf on $\PP^n=\PP^n_\KK$.
The \emph{cohomology table} of $\cF$ is the collection of numbers
$$ 
\gamma(\cF) =(\gamma_{i,d}) \hbox{ with } \gamma_{i,d}=\dim H^i(\PP^n,\cF(d)),
$$
which we think of as an element of the real vector space 
$\prod_{d=-\infty}^\infty \RR^{n+1}$.

In Eisenbud-Schreyer \cite{ES-BNC} we characterized the cohomology tables of vector bundles on 
$\PP^n$ (up to a positive rational multiple) as the finite positive rational linear combinations of cohomology tables of
\emph{supernatural} bundles, which we described explicitly. In this paper we treat the 
cohomology tables of all coherent sheaves. These are given by infinite sums:

\begin{theorem}
The cohomology table of any coherent sheaf on $\PP^n$
can be written as a convergent series,
with positive real coefficients, of cohomology
tables of supernatural bundles supported on linear subspaces.
\end{theorem}

We actually prove a more precise result, which includes a uniqueness statement. To state
it we recall some ideas from Eisenbud-Schreyer \cite{ES-BNC}.

A sheaf $\cF$ on $\PP^{n}$ has \emph{supernatural cohomology}
if, for each integer $d$, the cohomology $H^i(\cF(d))$ is nonzero for at most
one value of $i$ and, in addition,
the Hilbert polynomial $d\mapsto \chi(\cF(d))$
has distinct integral roots. We define the
\emph{root sequence} of a supernatural sheaf $\cF$ to be the 
sequence of roots of the Hilbert polynomial,
written in decreasing order, $z_1>\cdots> z_{s}$ where $s$ is the 
dimension of the support of $\cF$.
It will be convenient to put $z_0=\infty$ and $z_{s+1}=z_{s+2}\ldots =-\infty$.

The Hilbert polynomial and the cohomology table  of a supernatural sheaf $\cF$ 
are determined by the root sequence $(z_1,\dots,z_s)$ and the degree of $\cF$ as follows.
It is immediate that
$$ 
\chi(\cF(d))= \frac{\deg \cF}{s!} \prod_{i=1}^s (d-z_i).
$$
By Theorem 6.4 of our \cite{ES-BNC},
 $$
 h^j\cF(d)=
 \begin{cases}
  \frac{\deg \cF}{s!} \prod_{i=1}^{s}\mid d-z_i\mid &\text{if $z_{j}>d> z_{j+1}$},\\
 0 & \text{otherwise}.
 \end{cases}
$$
By  Theorem  6.1 of that paper, there exists a supernatural sheaf
of dimension $s$ and degree $s!$ with any given root
sequence $z=(z_1>\cdots>z_s)$. It may be taken to be a vector
bundle on a linear subspace $\PP^s\subset \PP^n$.
We denote its cohomology table by $\gamma^z$.
Thus 
$\gamma^z$ is the cohomology table of a vector bundle if and only if 
$z_n>-\infty$.

We partially order the root sequences termwise, setting $z \ge z'$ when
$$
z_1 \ge z_1', \ldots,  z_n \ge z'_n.
$$
By a \emph{chain} we mean a totally ordered set.
If $Z$ is an infinite sequence of root sequences, $(q_z)_{z\in Z}$ a sequence of numbers, and $\gamma$ is a cohomology table, we write
$
\gamma= \sum_{z\in Z} q_z \gamma^z,
$
to mean
that each entry $\sum_{z\in Z} q_z \gamma^z_{i,d}$ converges to $\gamma_{i,d}$.

With these preparations we can state the precise version of our main result. 
Recall that a sheaf is said to be \emph{purely $s$-dimensional} if all its associated
subvarieties have dimension exactly $s$.

\begin{theorem}\label{main}
Let $\gamma(\cF)$ be the cohomology table of a coherent sheaf $\cF$ on $\PP^n$.
There is a chain of zero-sequences $Z$ and positive real numbers $q_z$ such that
$$
\gamma(\cF) = \sum_{z\in Z} q_z \gamma^z.
$$
Both $Z$ and the numbers $q_z$ are uniquely determined by this condition.
The coefficients $q_z$ corresponding 
 to cohomology tables $\gamma^z$ of dimension $\dim \cF$ are rational numbers.
 If $\cF$ is purely $s$-dimensional,
 then  all the $\gamma^z$ are cohomology tables of vector bundles on $\PP^s$
 and all the $q_z$ are rational.
 \end{theorem}

We do not know whether all the numbers $q_z$ are rational, nor whether,
if all the $\gamma^z$ are cohomology tables of vector bundles,
the sheaf $\cF$ is necessarily torsion-free.

When we want to display (parts) of a cohomology table we use the convention
\small $$
\begin{matrix}
\cdots&\gamma_{n,-n-1}&\gamma_{n,-n}&\gamma_{n,-n+1}&\cdots&\vline&n\cr
&\vdots&\vdots&\vdots &&\vline&\vdots\cr
\cdots&\gamma_{1,-2}&\gamma_{1,-1}&\gamma_{1,0}&\cdots&\vline&1\cr
\cdots&\gamma_{0,-1}&\gamma_{0,0}&\gamma_{0,1}&\cdots&\vline&0& \cr
\hline
\cdots&-1&0&1&\cdots&\vline&d\backslash i \cr
\end{matrix}
$$
\normalsize
We make this choice of indexing so that the cohomology table
of a coherent sheaf $\cF$ coincides with
the Betti table of the \emph{Tate resolution} of $\cF$. 
This is a minimal, doubly
infinite, exact free complex over the exterior algebra on $n+1$ generators that is
associated to $\cF$ by the Bernstein-Gel'fand-Gel'fand correspondence. It is studied in
Eisenbud-Fl\o ystad-Schreyer \cite{EFS} and Eisenbud-Schreyer \cite{ES-C}.
For consistency with the notation of those papers,
we number the rows from the bottom and the columns from left to right
as in the table above.

\begin{example} The ideal sheaf $\cI_p$ of a point in $\PP^2$
has the cohomology table
\small
\setcounter{MaxMatrixCols}{20}
$$
\begin{matrix}
\cdots&10&6&3&1&  &  &  &   &          &&\vline&2\cr
\cdots&1&1&1&1&1&  &  &  &           &&\vline&1\cr
           & &   &   &  &  &2&5&9&14&\cdots&\vline&0& \cr
           \hline
\cdots&-4&-3&-2&-1&0&1&2&3&4&\cdots&\vline&d\backslash i \cr
\end{matrix}
$$
\normalsize

\noindent
where we drop the zero entries to make the shape more visible.
The expression in Theorem \ref{main} is
$$
\gamma(\cI_p)= \sum_{k=2}^\infty q_{(0,-k)}\gamma^{(0,-k)}
$$
where 
$$
q_{(0,-k)}=\frac{2}{(k-1)k(k+1)}.
$$
In particular
$$
\sum_{k=2}^\infty \frac{2d(d+k)}{(k-1)k(k+1)} = {d+2 \choose 2} -1
$$
holds for any $d \ge 1$,
$$
\sum_{k=-d+1}^\infty \frac{2d(d+k)}{(k-1)k(k+1)}=-1
$$
for any $d \le-1$
and 
$$
\sum_{k=2}^{-d} \frac{2d(d+k)}{(k-1)k(k+1)}=\frac{(d+2)(d+1)}{2}
$$
for any $d \le-2$.
\end{example}

To explain the proof of Theorem \ref{main}, we introduce a little more terminology.
We define the  $i$-th  {\it regularity} of a table $\gamma\in \prod_{d=-\infty}^\infty \RR^{n+1}$
to be 
$$
z_i(\gamma)=\inf \{ d \mid \gamma_{j,e+j} = 0 \hbox{ for all } j \ge i, \ e\geq d\}.
$$
We refer to 
$
z(\gamma) =(z_1(\gamma),\ldots,z_n(\gamma))
$ 
as the {\it regularity sequence} of  $\gamma$. It follows immediately from the definition that $z_1(\gamma) > z_2(\gamma) >\cdots$. 
Note that  $z_1(\gamma(\cF))$ coincides with the Castelnuovo-Mumford regularity of the sheaf $\cF$.
If $\gamma$ is the cohomology table of a supernatural sheaf $\cF$, then it follows from
Theorem 6.4 of our \cite{ES-BNC} that
$z_i(\gamma)$ is the $i$-th root of the Hilbert polynomial of $\cF$.

We define the \emph{support} of a table $\gamma$ to be the set of indices
$\{(i,d) \mid \gamma_{i,d}\neq 0\}$, and 
 the \emph{dimension} of $\gamma$ to be the
maximum $i$ such that $\gamma_{i,d}\not=0$ for some $d$, or $-1$ if 
all the $\gamma_{i,d}$ are zero.
Finally, the  \emph{corners} and \emph{corner values} of $\gamma$ are defined to be the positions
$$
(i,z_i(\gamma)+i-1) \hbox{ and values } \gamma_{i,z_i(\gamma)+i-1}
$$
for each $i$ such that  $i \le \dim \gamma$ and $z_{i+1} < z_i-1$.
 The decomposition of Theorem \ref{main} is effected by a transfinite ``greedy algorithm'': 

\begin{algorithm}{\bf (Decompose a Cohomology Table)}\label{Decomposition algorithm}\\
{\it Input:} A cohomology table $\gamma = \gamma(\cF)$ for some 
coherent sheaf $\cF$ on $\PP^n$.\\
{\it Output:} A chain of root sequences $Z$ and positive 
real numbers $(q_z)_{z\in Z}$
such that $\gamma = \sum_{z\in Z} q_z \gamma^z$.
\begin{enumerate}
\item Set $Z=\{\}$.
\item Set $i= \dim \gamma$.
\item WHILE $\dim \gamma = i$ DO
\begin{enumerate}
\item Let $z$ be the regularity sequence of $\gamma$, and replace $Z $ by $Z \cup \{z\}$
\item Let  $q_z>0$ be largest 
real number  such that the corner values of $\gamma$ are $\ge$ to the corner values of $q_z\gamma^z$. 
\item Replace $\gamma$ by $\gamma-q_z\gamma^z$.
\end{enumerate}
\item Replace $\gamma$ by the limit of the tables produced in step 3c.
\item If $\gamma=0$ then STOP, else go to Step 2.
\end{enumerate}
\end{algorithm}

Note that Step 2 is executed at most $n$ times, but we may loop through Steps 3a--3c 
infinitely often for each value of $i$ from $n$ to 1.
\smallskip

\begin{proof}[Outline of the proof that Algorithm \ref{Decomposition algorithm} succeeds] 
The crucial difficulty in the proof
of Theorem \ref{main}
is to show that  table $\gamma-q_z\gamma^z$ produced each time we pass through 
Step 3c has non-negative entries, and is sufficiently ``like'' the cohomology table of a coherent
sheaf to allow us to continue. To do this we will define a class of tables
closed under the basic operation in Step 3, and under taking limits in an appropriate way.
We call these \emph{admissible} tables; they are defined in \S \ref{Subtracting Once}.

The proof that Step 3c produces an admissible table is also given in \S \ref{Subtracting Once}.
It rests on an understanding of some functionals that are positive
on the cohomology tables of sheaves. 
Some of these functionals were defined in our paper \cite{ES-BNC},
and \S \ref{positivity} contains a simplified description of them, as well as some
others necessary for the present proof. 

The dimension $s$ of $\gamma$ is genuinely reduced each time we return to Step 2:
Indeed, some corner value of $\gamma$ becomes zero in Step 3c, decreasing some $z_i$.
Since $z_s$ remains the smallest of the (finite)
$z_i$, only finitely many steps can occur before $z_s$ is reduced, and thus in the course of
the WHILE loop, $z_s$ must be reduced to $-\infty$, so the dimension drops in Step 4, if it has
not dropped already in Step 3.

The convergence of the limiting process in Step 4 is dealt with in \S \ref{main proof},
as are the uniqueness and the special case of a pure-dimensional sheaf. 
Finally, the necessary positivity  is proven in
\S \ref{positivity section}, following an idea suggested by Rob Lazarsfeld. 
\end{proof}

The following example shows that the decomposition of 
Theorem \ref{main} sometimes mixes the torsion and torsion-free parts of a sheaf,
even when the sheaf itself is a direct sum.

\begin{example} Let $\cI$ be the ideal sheaf of a point in $\PP^2$,
and let $L$ be a line in $\PP^2$. Set $\cF = \cI \oplus \cO_L(-4)$.
The cohomology table of $\cF$ is given by the following diagram, where we have marked the corner values
with boxes.
$$
\begin{matrix}
\cdots&6&3&{\rlap{\hbox to 10pt{\hfil\raise0.1pt\hbox{$1$}\hfil}}}{\hbox to 10pt{\hskip-.9pt\vrule\hss\vbox to 10pt{\vskip-.2pt
            \hrule width 10pt\vfill\hrule\vskip-.8pt}\hss\vrule\hskip-.2pt}} 
&  &  &  &   &          &&\vline&2\cr
\cdots&8&7&6& 5 &3&2& {\rlap{\hbox to 10pt{\hfil\raise0.5pt\hbox{$1$}\hfil}}}{\hbox to 10pt{\hskip-.8pt\vrule\hss\vbox to 10pt{\vskip-.2pt
            \hrule width 10pt\vfill\hrule\vskip-.8pt}\hss\vrule\hskip-.4pt}}       &          &&\vline&1\cr
           &   &   &  &  &2&5&9&15&\cdots&\vline&0& \cr
           \hline
\cdots&-3&-2&-1&0&1&2&3&4&\cdots&\vline&d\backslash i \cr
\end{matrix}
$$
The regularity sequence is $z=(-2,3)$.
The supernatural cohomology table $\gamma^z$ is
$$
\begin{matrix}
\cdots&24&14&{\rlap{\hbox to 10pt{\hfil\raise0.1pt\hbox{$6$}\hfil}}}{\hbox to 10pt{\hskip-.9pt\vrule\hss\vbox to 10pt{\vskip-.2pt
            \hrule width 10pt\vfill\hrule\vskip-.8pt}\hss\vrule\hskip-.2pt}} 
&  &  &  &   &          &&\vline&2\cr
&&&& 4 &6&6& 
{\rlap{\hbox to 10pt{\hfil\raise0.5pt\hbox{$4$}\hfil}}}{\hbox to 10pt{\hskip-.8pt\vrule\hss\vbox to 10pt{\vskip-.2pt
            \hrule width 10pt\vfill\hrule\vskip-.8pt}\hss\vrule\hskip-.4pt}}       &          &&\vline&1\cr
           &   &   &  &  &&&&6&\cdots&\vline&0& \cr
           \hline
\cdots&-3&-2&-1&0&1&2&3&4&\cdots&\vline&d\backslash i \cr
\end{matrix}
$$
so we see that $q_z=1/6$. The table $\gamma':=\gamma-q_z\gamma^z$ has the form
$$
\begin{matrix}
\cdots&2&
{\rlap{\hbox to 20pt{\hfil\raise0.8pt\hbox{$2/3$}\hfil}}}{\hbox to 20pt{\hskip-.8pt\vrule\hss\vbox to 11pt{\vskip-.5pt
            \hrule width 20pt\vfill\hrule\vskip-.8pt}\hss\vrule\hskip-.4pt}} 
& &  &  &  &   &          &&\vline&2\cr
\cdots&8&7&6& 13/3 &2&1& {\rlap{\hbox to 20pt{\hfil\raise0.5pt\hbox{$1/3$}\hfil}}}{\hbox to 20pt{\hskip-.8pt\vrule\hss\vbox to 11pt{\vskip-.2pt
            \hrule width 20pt\vfill\hrule\vskip-.8pt}\hss\vrule\hskip-.4pt}}       &          &&\vline&1\cr
           &   &   &  &  &2&5&9&14&\cdots&\vline&0& \cr
           \hline
\cdots&-3&-2&-1&0&1&2&3&4&\cdots&\vline&d\backslash i \cr
\end{matrix}
$$
The regularity sequence of this table is $z'=(-3,3)$. This time, the corner that is cancelled
in $\gamma'$ is the one in the middle row, which comes from the torsion sheaf $\cO_L(-4)$,
rather than from $\cI$, and the table $\gamma'-q_{z'}\gamma^{z'}$ looks like
$$
\begin{matrix}
\cdots&14/15&{\rlap{\hbox to 20pt{\hfil\raise0.1pt\hbox{$1/5$}\hfil}}}{\hbox to 20pt{\hskip-.9pt\vrule\hss\vbox to 11pt{\vskip-.2pt
            \hrule width 20pt\vfill\hrule\vskip-.8pt}\hss\vrule\hskip-.2pt}} 
& &  &  &  &   &          &&\vline&2\cr
\cdots&8&7&17/3& 19/5 &7/5& {\rlap{\hbox to 25pt{\hfil\raise0.5pt\hbox{$7/15$}\hfil}}}{\hbox to 25pt{\hskip-.8pt\vrule\hss\vbox to 11pt{\vskip-.2pt
            \hrule width 25pt\vfill\hrule\vskip-.8pt}\hss\vrule\hskip-.4pt}}        &&         &&\vline&1\cr
           &   &   &  &  &2&5&9&203/15&\cdots&\vline&0& \cr
           \hline
\cdots&-3&-2&-1&0&1&2&3&4&\cdots&\vline&d\backslash i \cr
\end{matrix}
$$
\end{example}

\noindent{\bf Acknowledgements:}
We have enjoyed discussions of the material here with
Mats Boij and Rob Lazarsfeld.
We are particularly grateful to Lazarsfeld, who opened the path to this paper by pointing
out that one could use a \v Cech complex instead of a monad in the proof of
Theorem \ref{positive} for vector bundles. The program Macaulay2 \cite{M2}
of Mike Stillman and Dan Grayson has, once again, been invaluable in collecting evidence for
our conjectures and in suggesting how the proofs might go.  Finally, Silvio Levy helped us, with his
usual generosity, with advice on exposition and expertise about TeX. 

\section{Positive Functionals on Cohomology Tables}
\label{positivity}

In this section we will define some functionals---that is, real valued functions---
of an array
$$
\gamma = (\gamma_{j,d}) \in 
\prod_{d=-\infty}^\infty \RR^{n+1}.
$$
The key to the proof of the Theorem \ref{main} is the
Positivity Theorem \ref{positive} below, stating that
certain of these functionals take non-negative values on the 
cohomology tables of coherent sheaves.

Some of the functionals we need were defined in our \cite{ES-BNC}, and
Theorem \ref{positive} for those functionals, in the case of the cohomology table of a vector
bundle, is a translation of what is there.
Here we present a much simpler account of the functionals, that adapts well to
the new ones we use. The proof of Theorem \ref{positive}  given in \S \ref{positivity section}.

Define the $t$-th \emph{partial Euler characteristic of the $d$-th twist}
 of a table $\gamma\in \prod_{d=-\infty}^\infty \RR^{n+1}$  to
be the functional
$$
\chi_d^{\leq t}(\gamma) = \sum_{i=0}^t (-1)^i\gamma_{i,d}.
$$
When $t=\infty$ (or is simply large enough to be irrelevant) we
simply write  $\chi_d(\gamma)$ instead of $\chi_d^{\leq t}$.
For example, the usual Euler characteristic of a sheaf $\cF$ on $\PP^n$
is $\chi(\cF)=\chi_0^{\leq n}(\gamma(\cF))=\chi_0(\gamma(\cF)) $.

If 
\begin{eqnarray*}
d &=&d_0, \dots, d_{s+1} \in \ZZ,\cr
\psi &=&\psi_0,\dots, \psi_{s+1}\in \ZZ\cup\{\infty\}
\end{eqnarray*}
are sequences 
 (which we will call \emph{degrees} and \emph{bounds}, respectively)
we set 
$$
r_i=r_i(d):= \prod_{\substack{ 0\le j < k \le s+1\\j,k \neq i }}(d_k-d_j)
$$
and define a functional
$$
L(d,\psi) \colon  \prod_{d \in \ZZ} \RR^{n+1} \to \RR
$$ 
by the formula
\begin{align*}
L(d,\psi)&=\sum_{i=0}^{s+1} (-1)^ir_i\chi_{-d_i}^{\leq \psi_i}\\
\gamma&\mapsto
\sum_{i=0}^{s+1} (-1)^ir_i\chi_{-d_i}^{\leq \psi_i}(\gamma) = \sum_{i=0}^{s+1} (-1)^i r_i \sum_{j=0}^{\psi_i}(-1)^j\gamma_{j, - d_i} .
\end{align*}
We write $\infty$ for the special sequence of bounds $(\infty,\dots,\infty)$.
The naturalness of the functionals $L(d,\psi)$ is suggested by the following
well-known result used for interpolating polynomials, and its specialization to our case.

\begin{lemma}\label{interpolation} Let $d=(d_0,\dots,d_{s+1})$ be any sequence
of $s+2$ numbers, and let $r_i=r_i(d)$ as above. 
If
 $\gamma$ is the cohomology table of a coherent sheaf of dimension $\leq s$
 (or any table of dimension $s$
  such that $d\mapsto \chi_d(\gamma)$ is a polynomial of degree $\leq s$)
 then 
$
L(d,\infty)(\gamma) = 0.
$
\end{lemma}

\begin{proof} More generally, 
if $p(t)$ is
any polynomial of degree $\leq s$, then
$$
\sum_{i=0}^{s +1}(-1)^i r_i p(d_i) = 0.
$$
This follows from the fact that the last column of the $(s+2)\times(s+2)$ matrix
\scriptsize
$$
\begin{pmatrix}
1&d_0&\cdots &d_{0}^s&p(d_{0})\cr
1&d_1&\cdots&d_{1}^s&p(d_{1})\cr
\vdots&\vdots&&\vdots&\vdots\cr
1&d_s&\cdots&d_{s}^s&p(d_{s}) \cr
1&d_{s+1}&\cdots &d_{s+1}^s&p(d_{s+1}) \cr
\end{pmatrix}
$$
\normalsize
is linearly dependent on the others, so the determinant vanishes.
The displayed formula is the Laplace expansion of this determinant along the last
column.
\end{proof}

We will use the $L(d,\psi)$ with some other special sequences of bounds $\psi=\phi^j$ as well.
They are defined as follows:
For $j=1,\dots,s$, we define
\begin{eqnarray*}
\phi^j(s)&=&(\phi^j_0,\dots \phi^j_{s+1}),\cr\cr
\hbox{where}&&\cr
\phi^{j} _i &=&
\begin{cases} 
i & \hbox{ if } i < j \cr
i-1 & \hbox{ if } i = j \cr
i-2 & \hbox{ if } i >j ,
\end{cases}
\end{eqnarray*}
or, less formally,
$$
\phi^j(s) = (0,\dots,j-2, j-1, j-1, j-1, j,\dots, s-1).
$$
Finally,  we set
$\phi^0(s) = (-1, 0,\dots, s-2, s-1, s-1)$.
Here is our main result on the functionals $L(d,\phi^j(s))$:

\begin{theorem}[Positivity]\label{positive}
Let $d$ be a degree sequence, $d=(d_0<\cdots<d_{s+1})$ and let $r=r(d)$.
If $\cF$ is a coherent sheaf on $\PP^n$, then, for all $j\geq 1$
$$
L(d,\phi^j(s))(\gamma(\cF)) \ge 0,
$$
and
$$
-L(d,\phi^0(s))(\gamma(\cF)) \ge 0.
$$
\end{theorem}

\begin{figure}[h]
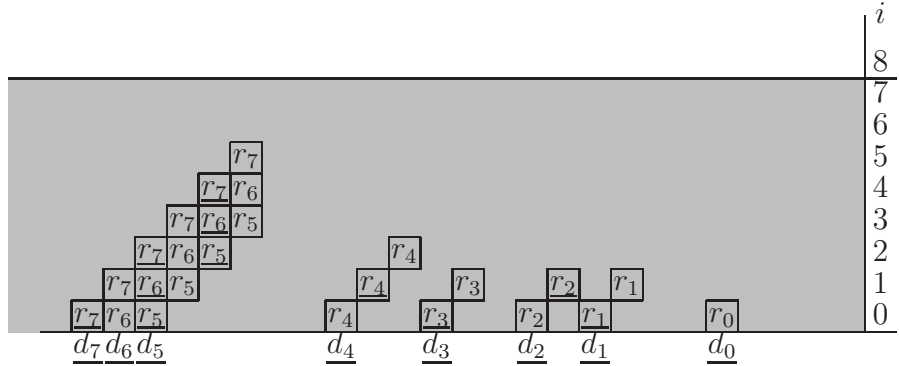

\begingroup \color{Gray} 
\leavevmode
\backgroundo
\jscaleo
\Plop $\r_7$ at (1,0)
\Plop $r_6$ at (2,0)
\Plop $r_7$ at (2,1)
\Plop $\r_5$ at (3,0)
\Plop $\r_6$ at (3,1)
\Plop $\r_7$ at (3,2)
\Plop $r_5$ at (4,1)
\Plop $r_6$ at (4,2)
\Plop $r_7$ at (4,3)
\Plop $\r_5$ at (5,2)
\Plop $\r_6$ at (5,3)
\Plop $\r_7$ at (5,4)
\Plop $r_5$ at (6,3)
\Plop $r_6$ at (6,4)
\Plop $r_7$ at (6,5)
\Plop $r_4$ at (9,0)
\Plop $\r_4$ at (10,1)
\Plop $r_4$ at (11,2)
\Plop $\r_3$ at (12,0)
\Plop $r_3$ at (13,1)
\Plop $r_2$ at (15,0)
\Plop $\r_1$ at (17,0)
\Plop $\r_2$ at (16,1)
\Plop $r_1$ at (18,1)
\Plop $r_0$ at (21,0)
\plop $\d_7$ at (1,-1)
\plop $\d_6$ at (2,-1)
\plop $\d_5$ at (3,-1)
\plop $\d_4$ at (9,-1)
\plop $\d_3$ at (12,-1)
\plop $\d_2$ at (15,-1)
\plop $\d_1$ at (17,-1)
\plop $\d_0$ at (21,-1)
\caption{To save space we have denoted $-d_i$ by $\d_i$ and  $-r_i$ by $\r_i$.
The shaded space indicates the positions where a cohomology table of a sheaf of dimension $7$ on 
$\PP^8$ could
have nonzero values. The functional $L(d,\phi^2(6))$ is the dot product with the table having
$\pm r_i$ in the positions shown, which are
initial segments of the diagonals numbered $\d_0,\ldots,\d_{7}$, and zeros elsewhere.
}
\label{dot product table}
\end{figure}

One may visualize the action of the linear form $L(d,\phi^j(s))$
 on a cohomology table $\gamma$
as the dot product of $\gamma$ with the table illustrated (for the case
$s=6, j=2$)
 in Figure \ref{dot product table}. 

For the case $j>0$
the proof, given in \S \ref{positivity section}, follows the same
outline as that in our paper \cite{ES-BNC}. 
Using the results of 
our \cite{ES-BNC} and Boij-S\"oderberg \cite{BS2}, 
Theorem \ref{positive}, in the case $j>0$, is
equivalent to Theorem \ref{coherent sheaf cohomology}.
We will deduce the case $j=0$ from the
case $j>0$ by a complicated numerical argument. It would be interesting to give 
a direct argument for the case $j=0$ as well.

Here is an example of how Theorem \ref{positivity} can be applied.
\begin{example} 
The Hilbert scheme $Hilb^{2t+2}(\PP^3)=H_1 \cup H_2$ has two irreducible components, which we
will call $H_1$ and $H_2$. 
The generic point of $H_1$ corresponds  to  two skew lines $X\subset \PP^3$, while
the generic point of $H_2$ corresponds to $Y=C\cup p\subset \PP^3$,
where $C$ is a conic and $p$ is a point not in the plane spanned by $C$.
The cohomology table of the ideal sheaf $\cI_X$ is
\setcounter{MaxMatrixCols}{20}
$$
\gamma(\cI_X)=\qquad
\begin{matrix}
\cdots&20&10&4&1&&    &  &   &          &&\vline&3\cr
\cdots&10&8&6&4&2&    &  &   &          &&\vline&2\cr
\cdots&  & & & & &1&   &  &           &&\vline&1\cr
           & &   &   &  &  &&4&12&25&\cdots&\vline&0& \cr
           \hline
\cdots&-4&-3&-2&-1&0&1&2&3&4&\cdots&\vline&d\backslash i \cr
\end{matrix}
$$
while that of $\cI_Y$ is
$$
\gamma(\cI_Y)=\qquad
\begin{matrix}
\cdots&20&10&4&1&&    &  &   &          &&\vline&3\cr
\cdots&11&9&7&5&3& 1   &  &   &          &&\vline&2\cr
\cdots&1 &1&1&1&1&1&   &  &           &&\vline&1\cr
           & &   &   &  &  &&4&12&25&\cdots&\vline&0& \cr
           \hline
\cdots&-4&-3&-2&-1&0&1&2&3&4&\cdots&\vline&d\backslash i \cr
\end{matrix}
$$
Using Theorem \ref{positivity} 
we can show that any integral table ``between" these two tables, obtained by
replacing the value $h^1\cI_Y(d)=1$ with a zero, and decreasing
$h^2\cI(d)$ by 1 as well, for some set of values $d<0$,
 cannot occur as the cohomology table of any sheaf; and
even that no multiple of such a table can occur. For example, no multiple of either the table
$$
T_2:=\qquad 
\begin{matrix}
\cdots&20&10&4&1&&    &  &   &          &&\vline&3\cr
\cdots&10&8&6&4&2& 1   &  &   &          &&\vline&2\cr
\cdots&  & & &&1&1&   &  &           &&\vline&1\cr
           & &   &   &  &  &&4&12&25&\cdots&\vline&0& \cr
           \hline
\cdots&-4&-3&-2&-1&0&1&2&3&4&\cdots&\vline&d\backslash i \cr
\end{matrix}
$$ 
or
$$
T_3:=\qquad 
\begin{matrix}
\cdots&20&10&4&1&&    &  &   &          &&\vline&3\cr
\cdots&10&8&6&4&3& 1   &  &   &          &&\vline&2\cr
\cdots&  & & &1&1&1&   &  &           &&\vline&1\cr
           & &   &   &  &  &&4&12&25&\cdots&\vline&0& \cr
           \hline
\cdots&-4&-3&-2&-1&0&1&2&3&4&\cdots&\vline&d\backslash i \cr
\end{matrix}
$$ 
can be the cohomology table of a coherent sheaf. 

One way to prove such a statement would be to apply Algorithm \ref{Decomposition algorithm},
and see that it eventually encounters a table with a negative entry.
For instance, in the case of the table $T_3$,
 that occurs after $16$  steps. But to prove the statement in general,
it is easier to appeal directly to Theorem \ref{positive}.

First, consider the functional
$L((-1,1,2,3),\phi^2(2))$, which may be written as the dot product with the table
$$
\begin{matrix}
&&&&&&    &  &   &          \vline&3\cr
&&&&&&    &  &   &          \vline&2\cr
&  & &6 &-16&12&&   &  &           \vline&1\cr
           & &  -6 & 16  & -12 &  &2&&&\vline&0& \cr
           \hline
\cdots&-4&-3&-2&-1&0&1&&\cdots&\vline&d\backslash i \cr
\end{matrix}
$$
in which all entries not shown are zero. The value of this functional on the table
$T_3$ shown above, for example, is $12-16 = -4$, proving that no multiple of
$T_3$ can be a cohomology table. Shifting this functional 
 to $L((-1+e, 1+e, 2+e, 3+e), \phi^2(2))$ we get a collection
 of functionals that prove the corresponding statement for any
 table between $\gamma(\cI_X)$ and $\gamma(\cI_Y)$ that has the pattern
 $0,1$  somewhere in the $h^1$ row, \emph{except} for $T_2$. However, the functional
$ L((-1,0,1,2,5),\phi^2(3))$, which is given by the dot product with the table
$$
\qquad 
\begin{matrix}
&&&&&&    &  &   &          \vline&3\cr
&&&12&&&    &  &   &          \vline&2\cr
&  &-12 & &&240&-540& 432  &  &   \vline&1\cr
            & 12  &   &  &-240  &540&-432&120&&\vline&0& \cr
           \hline
\cdots&-5&-4&-3&-2&-1&0&1&\cdots&\vline&d\backslash i \cr
\end{matrix}
$$ 
takes the value $432-540+12\cdot 8=-12$ on $T_2$, proving the claim.
\end{example}

\newpage
\section{Subtracting Once}
\label{Subtracting Once}


As we execute the the Algorithm \ref{Decomposition algorithm} we may leave the class
of cohomology tables of coherent sheaves. We will say that
a table is \emph{admissible} if it satisfies conditions 1-3 below.
As we shall see, the tables produced by the decomposition algorithm
will all be admissible.

The first two conditions that an admissible table $\gamma \in \prod \RR^{n+1}$ must
satisfy are: 
\begin{itemize}
\item[1.]
$\gamma_{i,d}=0$ for $i>0$ and $d\gg 0$. 
\item[2.]
The function 
$
d\mapsto \chi_d(\gamma)
$ 
from $\ZZ$ to $\RR$ is a polynomial of degree $s'\leq\dim \gamma$.
\end{itemize}
We will see that, in fact, admissibility implies that the degree of the polynomial
in condition 2 is exactly $\dim \gamma$ (Corollary \ref{equality of degrees}.)

For the last condition we need two definitions.
Suppose that $\gamma$ is a table satisfying 1 and 2. Suppose that
the dimension of $\gamma$ is $s$, and let $z_1>\cdots>z_s$ be the
regularity sequence of $\gamma$, as defined above. We call the 
table positions
$$
\{(i,d) \mid z_{i+1}< d+i<z_i\}
$$ 
the \emph{top positions} of $\gamma$, and all other positions with possibly nonzero values
$$
\{(i,d) \mid d+i \le z_{i+1}\}
$$
the \emph{lower positions} of $\gamma$. The last condition for admissibility is:
\begin{itemize}
\item[3.]
The values at the lower positions of $\gamma$ coincide with the values of the cohomology table of a coherent sheaf. That is, there exists a coherent sheaf $\cF$ such that
$$
\gamma_{i,d}= h^i(÷\cF(d)) \hbox{ for all lower positions } (i,d) \hbox{ of } \gamma
$$ 
\end{itemize}

\begin{figure}[sh]
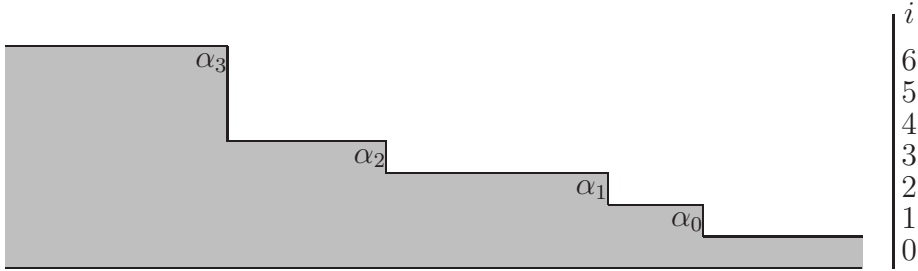

\begingroup \color{Gray} 
\leavevmode
\background
\jscale
\plop $\alpha_3$ at (6,6)
\plop $\alpha_2$ at (11,3)
\plop $\alpha_1$ at (18,2)
\plop $\alpha_0$ at (21,1)
\caption{A cohomology table of dimension $s=6$. 
The shaded region indicates where the table may have nonzero elements. 
The $\alpha_i$ are the corner values.}
\label{basic table}
\end{figure}

Now let $\gamma$ be an admissible table of  dimension $s$ with regularity sequence $z=z(\gamma)=(z_1,\ldots,z_s)$,
for example one whose shape is suggested by
Figure \ref{basic table}.

\begin{figure}[b]
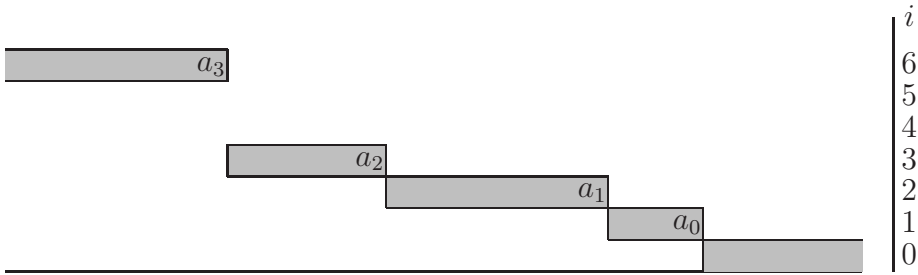

\leavevmode
\step 27 0 0    
\begingroup \color{Gray} 
\backtrack
\graystripe 7 7  
\graystripe 5 4
\graystripe 7 3
\graystripe 3 2
\graystripe 5 1
\endgroup
\backtrack
\step 7 7 6     
\step 5 4 3
\step 7 3 2
\step 3 2 1
\step 5 1 1     
\backtrack
\step 7 6 6     
\step 0 4 3
\step 5 3 2
\step 7 2 1
\step 3 1 0
\step 5 1 1
\backtrack
\jscale
\plop $a_3$ at (6,6)
\plop $a_2$ at (11,3)
\plop $a_1$ at (18,2)
\plop $a_0$ at (21,1)
\caption{Supernatural table $\gamma^z$ corresponding to the cohomology table in Figure \ref{basic table}. Here the $a_i$ are 
the corner values. The grayed area, where this table has nonzero values, coincides with the \emph{top positions} of the table in 
Figure \ref{basic table}.}
\label{basic supernatural table}
\end{figure}

We want to  subtract a suitable multiple $q_z\gamma^z$
of a supernatural table $\gamma^z$
so that, in $\gamma-q_z\gamma^z$, at least one of the corner values becomes zero, and the other
corner values
remain non-negative. Figures \ref{basic table} and \ref{basic supernatural table} give an idea of the pattern.

To achieve this goal we must take
$$
q_z=\min\{ \frac{\alpha_0}{a_0}, \dots, \frac{\alpha_m}{a_m}\},
$$
where $\alpha_0,\ldots,\alpha_m$ and $a_0,\ldots,a_m$ denote the corner values of $\gamma$ and $\gamma^z$ respectively.
The main step in the proof of Theorem \ref{main} is to show that \emph{all} of the entries of 
$\gamma-q_z\gamma^z$ are non-negative. This is the content of Proposition \ref{subtracting once}.

\begin{proposition}\label{subtracting once}
Let $\gamma$ be an admissible table of dimension $s>0$ with regularity sequence 
$
z=(z_1,\ldots,z_s)$. Let $\gamma^z$ be the cohomology table of a supernatural sheaf of dimension
$s=\dim \gamma$ with root sequence $z$. Let
$$
q_z=\min\{ \frac{\gamma_{i,z_i+i-1}}{\gamma^z_{i,z_i+i-1}} \mid i \le s \hbox{ and } z_{i+1} < z_i-1 \}
$$ 
be the minimal ratio of the corner values of $\gamma$ and $\gamma^z$. Then all entries of 
the table
$$
\gamma-q_z\gamma^z
$$ 
are non-negative, and its regularity sequence is $< z$.
\end{proposition}

\begin{corollary} \label{equality of degrees}
If $\gamma$ is a nonzero admissible table, then the function $d\mapsto \chi_d(\gamma)$
is a polynomial of degree exactly $\dim \gamma$.
\end{corollary}

\begin{proof}[Proof of Corollary \ref{equality of degrees}.]
For $d\gg 0$, the entry on the $d$-th diagonal of the table $\gamma^z$
is positive. Its value is 
$\prod_1^s(d-z_i)$,
and thus grows as a polynomial of degree $s=\dim \gamma^z = \dim \gamma$.
If $d\mapsto \chi_d(\gamma)$ had degree $<\dim \gamma$, then
$\gamma-q_z \gamma^z$ would have negative entries in these places,
contradicting Proposition \ref{subtracting once}.
\end{proof}

\begin{proof}[Proof of Proposition \ref{subtracting once}.]
Let $j$ be a cohomological index and $t$ a degree where
$\gamma^z_{j,t} \not=0$, say $z_{j+1}+j < t < z_j+j$. Let $\beta=\gamma_{j,d}$
and $b=\gamma^z_{j,d}$. We must show that 
 $\beta-q_z b \ge 0$.
 
If $t=z_j+j-1$ then we are talking about values at a corner position of $\gamma$ and $\gamma^z$,
and
the assertion follows immediately from the definition of $q_z$.
Thus we suppose that we are not at a corner position, that is,
$z_{j+1}+j< t < z_j+j-1$.  

We first treat the case where $j>0$. Figure \ref{case 1} illustrates the situation
for $j=2$.

\begin{figure}[h]
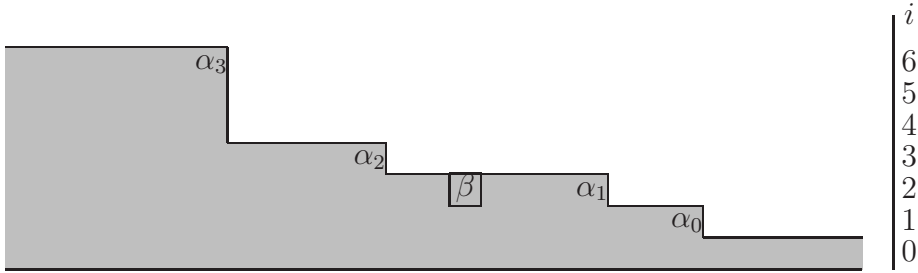

\begingroup \color{Gray} 
\leavevmode
\background
\jscale
\plop $\alpha_3$ at (6,6)
\plop $\alpha_2$ at (11,3)
\plop $\alpha_1$ at (18,2)
\plop $\alpha_0$ at (21,1)
\Plop $\beta$ at (14,2)
\caption{The case $j>0$ (here $j=2$). We must prove that the entry $\beta-q_zb$, of the table
$\gamma-q_z\gamma^z$, is non-negative.  Figure \ref{case 1 supernatural} shows the
corresponding entry of $\gamma^z$.} 
\label{case 1}
\end{figure}

\begin{figure}[b]
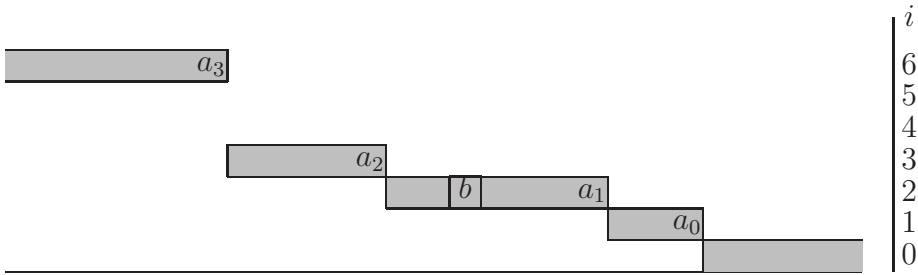

\leavevmode
\step 27 0 0    
\begingroup \color{Gray} 
\backtrack
\graystripe 7 7  
\graystripe 5 4
\graystripe 7 3
\graystripe 3 2
\graystripe 5 1
\endgroup
\backtrack
\step 7 7 6     
\step 5 4 3
\step 7 3 2
\step 3 2 1
\step 5 1 1     
\backtrack
\step 7 6 6     
\step 0 4 3
\step 5 3 2
\step 7 2 1
\step 3 1 0
\step 5 1 1
\backtrack
\jscale
\plop $a_3$ at (6,6)
\plop $a_2$ at (11,3)
\plop $a_1$ at (18,2)
\plop $a_0$ at (21,1)
\Plop $b$ at (14,2)
\caption{Supernatural table $\gamma^z$ showing the value $b$ at the same position as that of $\beta$ in
Figure \ref{case 1}.}
\label{case 1 supernatural}
\end{figure}

As indicated in the diagram, there is a corner
position of $\gamma$ and $\gamma_z$ immediately to the right of the position $(j,t)$,
and the values there are $\alpha_i:=\gamma_{j,z_j+j-1}$ and
$a_i:=\gamma^z_{j,z_j+j-1}$
respectively.
Since $\frac{\alpha_i}{a_i}\geq q_z$ it suffices to prove that
$$
\beta- \frac{\alpha_i}{a_i}b \geq 0 .
$$
To this end, consider the degree sequence
$$
d=(d_0,\ldots,d_{s+1}):=(-z_1,\ldots,-z_j,-z_j+1,-t+j ,-z_{j+1}, \ldots, -z_s)
$$
and let $r_i=r_i(d)$ as usual. Since $\chi_{z_i}(\gamma^z) = 0$ by construction,
Lemma \ref{interpolation} applied to the table $\gamma^z$ gives
$$
0=L(d,\infty)(\gamma^z) = \sum_{i=0}^{s +1}(-1)^i r_i \chi_{-d_i}(\gamma^z) = r_j a_i-r_{j+1}b, 
$$
so $b/a_i=r_j/r_{j+1}$, and it suffices to show that $r_{j+1}\beta - r_j\alpha_i \geq 0$.

On the other hand, we may apply Lemma \ref{interpolation} to the admissible table $\gamma$
to get
$$
0=L(d,\infty)(\gamma) = \sum_{i=0}^{s +1}(-1)^i r_i \chi_{-d_i}(\gamma) = r_j\alpha_i-r_{j+1}\beta+L(d,\phi^j)(\gamma).
$$
By the choice of the degree sequence $d$, the formula for $L(d,\phi^j)(\gamma)$ involves only
values at the lower positions of $\gamma$ (see Figure \ref{applied functional}.)

\begin{figure}[h]
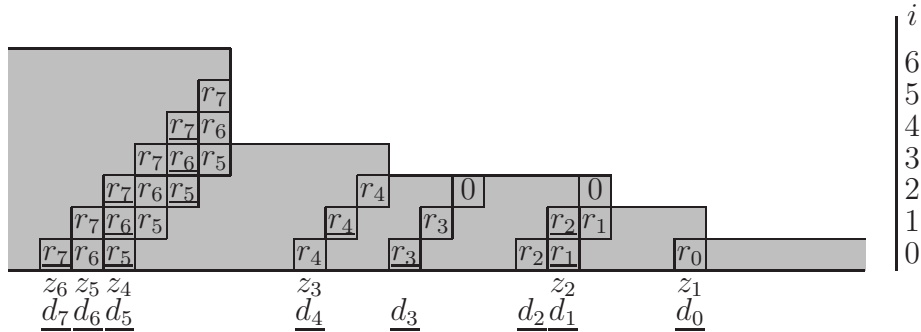

\begingroup \color{Gray} 
\leavevmode
\background
\jscale
\Plop $\r_7$ at (1,0)
\Plop $r_6$ at (2,0)
\Plop $r_7$ at (2,1)
\Plop $\r_5$ at (3,0)
\Plop $\r_6$ at (3,1)
\Plop $\r_7$ at (3,2)
\Plop $r_5$ at (4,1)
\Plop $r_6$ at (4,2)
\Plop $r_7$ at (4,3)
\Plop $\r_5$ at (5,2)
\Plop $\r_6$ at (5,3)
\Plop $\r_7$ at (5,4)
\Plop $r_5$ at (6,3)
\Plop $r_6$ at (6,4)
\Plop $r_7$ at (6,5)
\Plop $r_4$ at (9,0)
\Plop $\r_4$ at (10,1)
\Plop $r_4$ at (11,2)
\Plop $\r_3$ at (12,0)
\Plop $r_3$ at (13,1)
\Plop $$ at (14,2)
\lowerplop $0$ at (14,2)
\Plop $r_2$ at (16,0)
\Plop $\r_1$ at (17,0)
\Plop $\r_2$ at (17,1)
\Plop $r_1$ at (18,1)
\lowerplop $0$ at (18,2)
\Plop $$ at (18,2)
\Plop $r_0$ at (21,0)
\plop $z_6$ at (1,-1)
\plop $z_5$ at (2,-1)
\plop $z_4$ at (3,-1)
\plop $z_3$ at (9,-1)
\plop $z_2$ at (17,-1)
\plop $z_1$ at (21,-1)
\plop $\d_7$ at (1,-1.9)
\plop $\d_6$ at (2,-1.9)
\plop $\d_5$ at (3,-1.9)
\plop $\d_4$ at (9,-1.9)
\plop $\d_3$ at (12,-1.9)
\plop $\d_2$ at (16,-1.9)
\plop $\d_1$ at (17,-1.9)
\plop $\d_0$ at (21,-1.9)
\caption{The functional $L(d,\phi^j)$ is the dot product with the table having
$\pm r_i$ in the positions shown, and zeros elsewhere. In the illustration,  $s=6$ and $j=2$.
To save space we have denoted $-d_i$ by $\d_i$ and  $-r_i$ by $\r_i$. The explicit
zeros are added for emphasis.}
\label{applied functional}
\end{figure}

 Because $\gamma$ is admissible,
 $L(d,\phi^j)(\gamma)=L(d,\phi^j)(\gamma(\cF))$
for some coherent sheaf $\cF$. Thus we may apply Theorem \ref{positive} to conclude that
$$
r_{j+1}\beta - r_j\alpha_i =L(d,\phi^j)(\gamma) \ge 0
$$
as desired. \medskip 

The proof in the case $j=0$ is almost the same. Figure \ref{j=0 figure}
\begin{figure}[h]
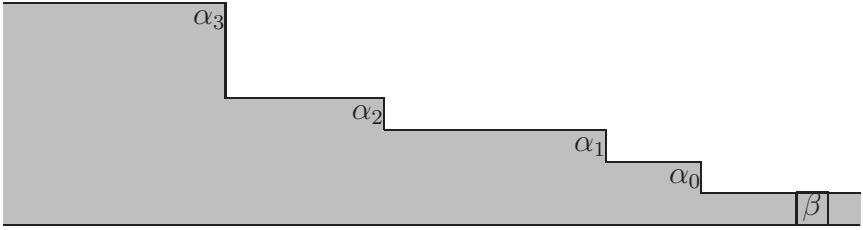

\begingroup \color{Gray} 
\leavevmode
\background
\plop $\alpha_3$ at (6,6)
\plop $\alpha_2$ at (11,3)
\plop $\alpha_1$ at (18,2)
\plop $\alpha_0$ at (21,1)
\Plop $\beta$ at (25,0)
\caption{Position of $\beta$ in case $j=0$.}
\label{j=0 figure}
\end{figure}
Illustrates the position of the value $\beta$ in this case.
Since $\gamma^z$ is
assumed nonzero at the position $(j,t)$, we must have 
 $t>z_1$ in this case. This time there is no corner position to the right of $(0,t)$,
 but we set $i=m$, and we let $d$ be the degree sequence
 $$
d=(d_0,\ldots,d_{s+1})=(-t,-z_1,\ldots,-z_s,-z_s+1).
$$
Figure \ref{applied functional 0} illustrates the relation of $\gamma$ to 
the positions involved in the functional $-L(d,\phi^0)$.
The rest of the argument is nearly the same. 
$$
\beta-q_zb \ge \beta - \frac{\alpha_m}{a_m}b  \ge 0
$$
follows, because
$$
0=L(d,\infty)(\gamma^z)=r_0b-r_{s+1}a_m
$$
gives ${b}/{a_m}={r_{s+1}}/{r_0}$, and
$$
0=L(d,\infty)(\gamma)=r_0\beta-r_{s+1}\alpha_m+L(d,\phi^0)(\gamma)
$$
implies the desired positivity, because $-L(d,\phi^0)(\gamma)\ge 0$ by Theorem \ref{positivity}.
\begin{figure}[h]
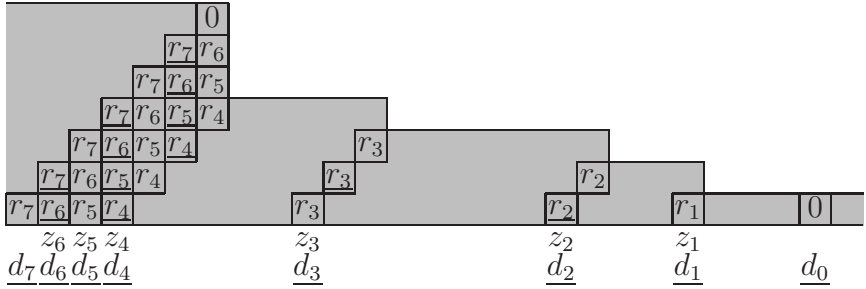

\begingroup \color{Gray} 
\leavevmode
\background
\step 27 0 0 
\backtrack
\Plop $r_7$ at (0,0)
\Plop $\r_6$ at (1,0)
\Plop $\r_7$ at (1,1)
\Plop $r_5$ at (2,0)
\Plop $r_6$ at (2,1)
\Plop $r_7$ at (2,2)
\Plop $\r_4$ at (3,0)
\Plop $\r_5$ at (3,1)
\Plop $\r_6$ at (3,2)
\Plop $\r_7$ at (3,3)
\Plop $r_4$ at (4,1)
\Plop $r_5$ at (4,2)
\Plop $r_6$ at (4,3)
\Plop $r_7$ at (4,4)
\Plop $\r_4$ at (5,2)
\Plop $\r_5$ at (5,3)
\Plop $\r_6$ at (5,4)
\Plop $\r_7$ at (5,5)
\Plop $r_4$ at (6,3)
\Plop $r_5$ at (6,4)
\Plop $r_6$ at (6,5)
\lowerplop $0$ at (6,6)
\Plop $$ at (6,6)
\Plop $r_3$ at (9,0)
\Plop $\r_3$ at (10,1)
\Plop $r_3$ at (11,2)
\Plop $\r_2$ at (17,0)
\Plop $r_2$ at (18,1)
\Plop $r_1$ at (21,0)
\lowerplop $0$ at (25,0)
\Plop $$ at (25,0)
%
\plop $z_6$ at (1,-1)
\plop $z_5$ at (2,-1)
\plop $z_4$ at (3,-1)
\plop $z_3$ at (9,-1)
\plop $z_2$ at (17,-1)
\plop $z_1$ at (21,-1)
\plop $\d_7$ at (0,-1.9)
\plop $\d_6$ at (1,-1.9)
\plop $\d_5$ at (2,-1.9)
\plop $\d_4$ at (3,-1.9)
\plop $\d_3$ at (9,-1.9)
\plop $\d_2$ at (17,-1.9)
\plop $\d_1$ at (21,-1.9)
\plop $\d_0$ at (25,-1.9)
\caption{The functional $-L(d,\phi^0)$ is the dot product with the table having
$\pm r_i$ in the positions shown, and zeros elsewhere. In the illustration,  $s=6$.
To save space we have denoted $-d_i$ by $\d_i$ and  $-r_i$ by $\r_i$. The explicit
zeros are added for emphasis.}
\label{applied functional 0}
\end{figure}

\end{proof}

\section{Proof of the main result}\label{main proof}

We start by describing the growth of dimensions of the cohomology groups $h^i(\cF(d)) $ for $d \gg 0$.

\begin{proposition}\label{growth}
Let $\cF$ be a coherent sheaf on $\PP^n$. 
For each $i=0,\dots, n$ there exists a polynomial $p^{i}_\cF \in \QQ[t]$
such that
$$
p^{i}_\cF(d) = h^i (\cF(d)) \hbox{ \quad for all\quad } d \ll 0.
$$
The degree of $p^i_\cF$ is $\leq i$, with equality
if and only if $\cF$ has an associated subvariety of dimension
$i$. In particular, if the dimension of the support of $\cF$ is $s$, then
$\deg p^s_\cF = s$. Furthermore, $\cF$ is pure-dimensional, if and only if 
$\deg p^i_\cF <s$ for every $i<s$.
\end{proposition}

\begin{proof} Let $M$ be a graded module over the polynomial ring 
$S=\KK[x_0,\dots,x_n]$ whose associated sheaf is $\cF$.
For $i>0$,
\begin{eqnarray*}
\oplus_d \Hom_\KK( H^i(\cF(d)),\KK)& =&\oplus_d \Ext^{n-i}(\cF(d),\omega_{\PP^n}) \cr
&= &\oplus_d \Ext^{n-i}(\cF,\cO(-n-1-d)) \cr
&=&\Ext^{n-i}_S(M,S(-n-1)). \cr
\end{eqnarray*}
Thus $p_{\cF}^i$ is the Hilbert polynomial of 
$\Ext^{n-i}_S(M, S(-n-1))$, 
so the degree of $p_{\cF}^i$ is one less than the Krull dimension of
$\Ext^{n-i}_S(M, S(-n-1)),
$
or,
equivalently, of $ \Ext^{n-i}_S(M, S)$.

The inequality $\deg p^i_\cF\leq \dim \cF$ now follows from the 
Auslander-Buchsbaum-Serre Theorem: after localizing $S$ at any prime $P$ of dimension $>i+1$
we get a regular local ring of dimension $<(n+1)-(i+1) = n-i$, so $\Ext^{n-i}(M,S)_P=0$. It follows
that $\dim \Ext^{n-i}(M,S)<i+1.$ Now suppose that $P$ is a prime of dimension exactly $i+1$.
By the Auslander-Buchsbaum formula, 
$P$ is associated to $M$ if and only if the projective dimension of $M_P$ is $i+1$, which is true
if and only if  $\Ext^{n-i}(M,S)_P\neq 0$. Since every associated prime of a graded
module is homogeneous, $P$ must correspond in this case to an associated subvariety
of $\cF$, proving the statement about equality. The rest of the 
Proposition follows.
\end{proof}

\begin{proof}[Proof of Theorem \ref{main}.]
For the first statement of the Theorem it suffices to show that Algorithm \ref{Decomposition algorithm}
succeeds. We have already seen in Proposition \ref{subtracting once} that Step 3c, starting
with an admissible table, always produces a new admissible table, and we have explained
in the Outline of the Proof in the Introduction why the dimension of $\gamma$ will
drop by at least 1 each time we reach Step 2. Thus it suffices to show that
if we start with an admissible table $\gamma$, then the sequence
of tables produced by the WHILE loop of Step 3 actually converges to an admissible
table, so that we can execute Step 4.

Convergence is no problem: By Proposition \ref{subtracting once},
 the tables stay admissible, and thus have only
non-negative terms throughout an instance of Step 3. Thus the values in a given position form
a decreasing, bounded below sequence.

To show that the limiting table produced in Step 4 is actually admissible,
suppose the rows of cohomological index $s'+1, \ldots, s$ are wiped out by
a pass through Step 3,
 while the $s'$-th row remains nonzero.  We have to show that the remaining
table $\gamma'$ is an admissible table of dimension $s'$. It is clear, in any case,
that $\gamma'$ satisfies Condition 1 of admissibility.

Since the rows with cohomological index $0, \ldots,s'$ survive, only finitely many corner values  with cohomological index $j\le s'$ are removed in the course of Step 3. So we may replace $\gamma$
with the admissible table that results from finitely many subtractions, and assume that no corner value with cohomological index $\leq s'$
becomes zero in the infinite sequence of subtractions leading to $\gamma'$. 
It follows that the values of in the lower positions of $\gamma'$ are the same
as those in the corresponding positions of $\gamma$; thus condition 3 of admissibility
is satisfied.
 
To complete the proof, we first note that the sequence of Hilbert functions
of the tables obtained by the successive subtractions converges decreasingly
to a function that takes non-negative real values at all $d\gg 0$. 
At every finite stage we subtract a polynomial of degree $s+1$,
so the $s+1$-st
difference function is zero. By continuity, it remains zero in the limit.
It follows that $d\mapsto \chi_d(\gamma')$,
is a polynomial function.
 
On the other hand, the values on the top row of $\gamma'$ at the positions $d\ll 0$ grow at most like a polynomial of degree
$s'$ since all values are bounded by the values of the corresponding row of $\gamma$.
The rows with cohomological degree $i < s'$ have for $d \ll 0$ the values of the original 
cohomology table of $\cF$.
By Proposition \ref{growth}, they grow with negative
$d$ as polynomials of degree $<s'$. Thus the Euler characteristic
$\chi_{d}(\gamma) $ is a polynomial in $d$ of degree $\leq s'$; that is, $\gamma'$ satisfies
condition 2 of admissibility.
This completes the proof that Algorithm \ref{Decomposition algorithm} succeeds, and
produces a decomposition of the desired kind.

To prove uniqueness, suppose that $Z$ and $W$ are both chains of root sequences,
and that 
$$
\gamma(\cF) = \sum_{z\in Z} q_z \gamma^z= \sum_{w\in W} r_w \gamma^w
$$
 with
$q_z$ and $r_w$ positive real numbers, where $Z$ is the chain produced by
Algorithm \ref{Decomposition algorithm}. Since $Z$, at least, is well ordered,
there is a largest element of $Z$ that does not appear in $W$, 
or appears with a different coefficient. We may as well
subtract the contributions of the terms corresponding to larger elements of $Z$,
which are the same for the two sums, and thus suppose that 
$$
\gamma = \sum_{z\in Z} q_z \gamma^z= \sum_{w\in W} r_w \gamma^w
$$
where $\gamma$ is an admissible table, and the maximal element $\overline z \in Z$
either does not appear in $W$, or appears with a different coefficient $r_{\overline z}\neq q_{\overline z}$.

Because the root sequence of $Z$ is the regularity sequence of $\gamma$, every
$w\in W$ must satisfy $w\leq\overline z$. If $\overline z$ itself is in $W$, but
$r_{\overline z}\neq q_{\overline z}$, then 
$\gamma - r_{\overline z} \gamma^{\overline z}$ has exactly the same corner positions
and regularity sequence as $\gamma$. But since $W$ is a chain, at least one of the 
corner positions of $\gamma$ is represented with the value zero in every one of the 
$\gamma^w$ for $\overline z \neq w \in W$, and we see that
$
\gamma - \sum_{w\in W} r_w \gamma^w \neq 0,
$
contradicting our hypothesis. 

Similarly, if $\overline z\notin W$ then, since there are only finitely many elements just
below $z$ in the poset of root sequences, there is some corner position
of $\gamma$ that is represented by the value zero in every $\gamma^w$ for $w\in W$,
so we can finish the argument in the same way. This proves uniqueness.

Note that the coefficients $q_z$ involved in any finite sequence of subtractions
 in Algorithm \ref{Decomposition algorithm}
 starting from a rational cohomology table are automatically rational.
This applies to all the $q_z$ corresponding to $\gamma^z$ of dimension $=\dim \cF$.
 
Now suppose that $\cF$ is a pure-dimensional sheaf.
It suffices to show that the decomposition is obtained as the
 limit of finite sequences of subtractions starting from the cohomology table
 of $\cF$ in this case.
 
Once again, let $\gamma'$ be the result of subtracting 
the cohomology tables of vector bundles on $\PP^s$, as in Algorithm \ref{Decomposition algorithm},
so that $s':=\dim \gamma'<\dim \gamma(\cF) = s$.
By Proposition \ref{growth}, the Euler characteristic of
the resulting table $\gamma'$ grows like a polynomial of degree $\leq s'-1$.
If $\gamma'$ were nonzero, we would get a contradiction to Corollary
\ref{equality of degrees}.  Thus 
$\gamma'=0$, completing the proof.
\end{proof}

\section{Proof of the Positivity Theorem}
\label{positivity section}

In our paper \cite{ES-BNC} we defined pairings
$$
\langle \beta, \gamma\rangle = 
\sum_{\{i,j,k\mid j\leq i\}}
(-1)^{i-j}\beta_{i,k}\gamma_{j,-k}.
$$ 
and 
\begin{eqnarray*}
\langle \beta, \gamma\rangle_{c,\tau} &=&
\sum_{\{i,j,k\mid j\leq i \hbox{\scriptsize\ and } (j<\tau \hbox{\scriptsize\ or } j \leq i-2)\}}
(-1)^{i-j}\beta_{i,k}\gamma_{j,-k}
\cr
&+&
\sum_{\{i,j,k,\epsilon \mid 0\leq \epsilon\leq 1,\ j=\tau,\ i=j+\epsilon,\ k\leq c+\epsilon \}}
(-1)^{i-j}\beta_{i,k}\gamma_{j,-k}.
\end{eqnarray*}
for
$\beta = (\beta_{i,k})\in \oplus_{-\infty}^\infty \RR^{n+2}$
and $\gamma \in \prod_{-\infty}^\infty \RR^{n+2}$, and 
$0\leq \tau \leq n, c\in \ZZ$. We showed that, if $\beta$
is the Betti table of a finitely generated graded module over
$S:=\KK[x_0,\dots,x_n]$ and $\gamma$ is the cohomology table
of a vector bundle $\cF$, or of a complex $E$ of 
free graded $S$-modules, supported in positive cohomological 
degrees, then 
$$
\langle \beta, \gamma\rangle\geq 0\quad \hbox{ and }\quad\langle \beta, \gamma\rangle_{c,\tau}\geq 0.
$$
Our proof for the vector bundle case reduced to the case of a free complex
by replacing the vector bundle with a free monad. Since the free monads
of coherent sheaves have terms in negative cohomological degrees,
this proof could not show that the pairing above was non-negative when $\cF$ is
a general coherent sheaf. After our paper was
finished, Rob Lazarsfeld pointed out to
us a variation on our proof in which the monad for $\cF$ is replaced
by an injective or flasque resolution of $\cF$. It turns out that, with one further idea, this idea yields
a proof of non-negativity that works for any coherent sheaf $\cF$.

\begin{theorem}\label{coherent sheaf cohomology}
Let $F$ be the minimal free resolution of a finitely generated graded
 $S$-module $M$. If  $\cF$ is a coherent sheaf on 
$\PP^{n}$, then
$$
\langle F,\cF \rangle \geq 0 \quad\hbox{and}\quad \langle F,\cF \rangle_{c,\tau} \geq 0
$$
\end{theorem}

\begin{proof}
The number $\langle F,\cF\rangle$ depends only
on the dimensions of the $H^j(\cF(-k))$ for $k\in \ZZ$, we may begin
by replacing $\cF$ with a ``general translate'' by an element of $PGL(n)$, to
make $\cF$ homologically transverse to the sheaf $\tilde M$, as proven by
Sierra \cite{Sierra} and by Miller and Speyer \cite{Miller-Speyer}. If we let
$G$ be a graded $S$-module such that
$\widetilde G = \cF$, this means that the modules $\Tor_i^S(M, G)$ have support only
at the irrelevant ideal for $i>0$.

Let $E: \oplus_\ell G[x^{-1}_\ell] \to \cdots$ be the \v Cech complex of $G$.
The homological transversality implies that the complex $F\otimes E^j$ has
homology only at $F_0\otimes E^j$, so the total complex of the double complex
$F\otimes E$ has homology only in non-negative cohomological degree. We can now 
proceed exactly as in the proofs of Theorems 3.1 and 4.1 of our \cite{ES-BNC}.
\end{proof}

We next describe a
simplification in the statement that makes use of
 the main results of our \cite{ES-BNC} and
of Boij-S\"oderberg \cite{BS2}, and also an extension of the statement that
will be crucial for the proof of Theorem \ref{main}.

Recall that a graded
Cohen-Macaulay $S$-module $M$ of codimension $s+1$
 is said to have a \emph{pure} resolution
with degree sequence $d=(d_0, \dots, d_{s+1})$ if the minimal
free resolution of $M$ has the form
$$
S(-d_0)^{r_0} \lTo S(-d_1)^{r_1}\lTo\cdots\lTo S(-d_{s+1})^{r_{s+1}}\lTo 0.
$$
In this case, $d_0<\cdots<d_{s+1}$, and  there is a positive rational
number $q$ such that each $r_i = q\cdot r_i(d)$, where, as in \S \ref{positivity section}
$$
r_i(d):= \prod_{\substack{ 1\le j < k \le s+1\\j,k \neq i }}(d_k-d_j)
$$
(See Herzog and K\"uhl \cite{HK}).

Together, our \cite{ES-BNC} and
Boij-S\"oderberg \cite{BS2} show that there is a graded Cohen-Macaulay $S$-module
with any given degree sequence $(d_0<\cdots<d_{s+1})$, and 
the Betti table of any graded $S$-module is a positive
rational linear combination of the Betti tables of Cohen-Macaulay modules with
pure resolutions. Thus to prove that the value of a bilinear functional such as those above
is non-negative, it suffices to treat the case where $\beta$ is the Betti table of 
a Cohen-Macaulay module with pure resolution, and if the resolution has degree sequence
$d$, one may as well assume
that $r_i=r_i(d)$ for every $i$ as well: that is, we may restrict our attention
to the functionals $\langle(\beta^d, \gamma\rangle_{c,\tau}$ with
$\beta^d$ to be the table with 
$$
\beta^d:\quad \beta_{i,j} =
\begin{cases} 
r_i(d) & \text{ if }j= d_i \text{ and } \\
0&\text{ otherwise.}
\end{cases}
$$
 For such $\beta^d$ we
 may re-write the definition given above in the form:
 \begin{align*}
 \langle\beta^d, \gamma\rangle_{c,\tau} = \sum_{i<\tau}&(-1)^i \  r_i(d)\  \chi_{-d_i}^{\leq i}(\gamma)\\
 +&(-1)^\tau \  r_\tau(d)\  \chi_{-d_\tau}^{\leq A}(\gamma)
 \\
 +&(-1)^{\tau+1} r_{\tau+1}(d)\ \chi_{-d_{\tau+1}}^{\leq B}(\gamma) 
 \\
 +\sum_{i>\tau+1}&(-1)^i \  r_i(d)\ \chi_{-d_i}^{\leq i-2}(\gamma)
\end{align*}
where
\begin{align*}
A&=\begin{cases}
\tau-1 &\text{ if } c<d_\tau\\
\tau & \text{ otherwise}
 \end{cases}\\
B &=  \begin{cases}
\tau-1 &\text{ if } c<d_{\tau+1}\\
\tau & \text{ otherwise}
 \end{cases}
\end{align*}
It follows that if $\tau \geq 1$ and $c<d_\tau$ then
$$
\langle\beta^d, \gamma\rangle_{c,\tau} = L(d, \phi^\tau)(\gamma)
$$
while if $c\geq d_{\tau+1}$ then 
$$
\langle\beta^d, \gamma\rangle_{c,\tau} = L(d, \phi^{\tau+1})(\gamma).
$$
Moreover, if the $\gamma_{i,j}$ are non-negative, as in 
any admissible table, and $d_\tau\leq c <d_{\tau+1}$ then, comparing
signs, we see that 
$$
\langle\beta^d, \gamma\rangle_{c,\tau} \geq
\langle\beta^d, \gamma\rangle_{d_{\tau}-1,\tau}=  L(d, \phi^\tau)(\gamma)
$$
so this case is not very useful. 

\begin{proof}[Proof of Theorem \ref{positive}]
The description above shows that the cases  $j>0$
follow from Theorem \ref{coherent sheaf cohomology}.
 
To simplify the notation for the case $j=0$ we 
 set $\psi = \phi^0 = (-1, 0,1,\ldots,s-2,s-1,s-1)$.
 We write 
$$
d^{(j)}=(d_1,\ldots,d_j) \qquad \hbox{ for } j=1,\ldots, s+1
$$
and
 $$
  \psi^{(j)} =\begin{cases}(0,1,\ldots,j-1) & \hbox{ for }j=1, \ldots, s ,\, \hbox{ and } \cr
  (0,1,\dots,s-1,s-1) &\hbox{ for }j=s+1 .
 \end{cases}
 $$
We will show that
\begin{equation}
 -L(d,\psi)= 
 \sum_{\ell=0}^{s+1} (-1)^{s-\ell}r_{s+1-\ell}(d)\chi_{-d_{s+1-\ell}}^{\leq \psi_{s+1-\ell}} = 
\sum_{k=0}^{s}  A_k \; L(d^{(s+1-k)},\psi^{(s+1-k)})
\end{equation}
where
$$
A_k=\prod_{1\le j \le s-k }(d_j-d_0) \prod_{\substack{ 1\le i < j \le s+1\\s+1-k<j }}(d_j-d_i).
$$
The coefficients $A_k$ are obviously non-negative. By Theorem \ref{coherent sheaf cohomology}, 
the forms
$L(d^{(s+1-k)},\psi^{(s+1-k)})$
take non-negative values on the cohomology tables of coherent
sheaves, so this will suffice to prove  Theorem \ref{positive}.

The  coefficient of $(-1)^{s-\ell}\chi_{-d_{s+1-\ell}}^{\leq \psi_{s+1-\ell}}$ on the right-hand side of Equation (1) is
$$
\sum_{k=0}^\ell \quad 
\bigg( \prod_{1\leq j\leq s-k} (d_j-d_0) \bigg) 
\bigg( \prod_{\substack{1\le i< j\le s+1\\ s+1-k< j}}(d_j-d_i) \bigg)
\bigg( \prod_{\substack{1\le i <j \le s+1-k\\i,j \not=s+1- \ell}}(d_j-d_i) \bigg).
$$
We will show that this  is $r_{s+1-\ell}(d)$.
The terms in the sum have a common factor (coming from the first and third factors in each term)
$$
\bigg( \prod_{1\leq j\leq s-\ell} (d_j-d_0) \bigg) 
\bigg( \prod_{\substack{1\le i <j \le s+1-\ell\\ i,j \not=s+1- \ell}} (d_j-d_i) \bigg)
 =
  \prod_{\substack{0\le i <j \le s+1-\ell\\i,j \not=s+1- \ell}} (d_j-d_i).
$$
After factoring this out, we get
$$
\sum_{k=0}^\ell \quad 
\bigg( \prod_{s-\ell+1\leq j\leq s-k} (d_j-d_0) \bigg) 
\bigg( \prod_{\substack{1\le i< j\le s+1\\ s+1-k< j}} (d_j-d_i) \bigg)
\bigg(\prod_{\substack{1\le i <j \le s+1\\i,j \not=s+1- \ell\\s+1-\ell<j\leq s+1-k}} (d_j-d_i)\bigg),
$$
which can be further factored as
$$
\bigg(\prod_{\substack{1\le i <j\leq s+1\\i,j \neq s+1- \ell\\s+1-\ell<j}} (d_j-d_i)\bigg)
\sum_{k=0}^\ell
\bigg( \prod_{s-\ell+1\leq j\leq s-k} (d_j-d_0)\bigg)
\bigg( \prod_{\substack{i=s+1-\ell\\s-k+1< j}} (d_j-d_i)\bigg).
$$
Applying the case $t=-1$ of Lemma \ref{messy}, we can combine all the factors
to express the original sum as
$$
\bigg( \prod_{\substack{0\le i <j \le s+1-\ell\\i,j \not=s+1- \ell}} (d_j-d_i)\bigg)
\bigg(\prod_{\substack{1\le i <j\leq s+1\\i,j \neq s+1- \ell\\s+1-\ell<j}} (d_j-d_i)\bigg)
\bigg( \prod_{s-\ell+1< j\leq s+1} (d_j-d_0)\bigg)
$$
$$
=\prod_{\substack{0\le i<j\le s+1\\ i,j \not= s+1-\ell}}(d_j-d_i) = r_{s+1-\ell}(d),
$$
completing the proof.
\end{proof}

\begin{lemma}\label{messy}
For $-1\leq t\leq \ell-1$ we have
$$
\sum_{k=0}^\ell
\bigg( \prod_{s-\ell+1\leq j\leq s-k} (d_j-d_0)\bigg)
\bigg( \prod_{\substack{s-k+1< j}} (d_j-d_{s-\ell+1})\bigg)
$$
$$
=
\bigg( \prod_{s-\ell+1< j\leq s-t} (d_j-d_0)\bigg)
\bigg( \prod_{\substack{s-t< j}} (d_j-d_{s-\ell+1})\bigg)
$$
$$+
\sum_{k=0}^t
\bigg( \prod_{s-\ell+1\leq j\leq s-k} (d_j-d_0)\bigg)
\bigg( \prod_{\substack{s-k+1< j}} (d_j-d_{s-\ell+1})\bigg).
$$
\end{lemma}

\begin{proof}
The formula is obvious for $t=\ell-1$, so we do descending induction.
The induction step follows by combining the first product with the $k=t$ term
of the summation, as follows:
$$
\prod_{j=s-\ell+2}^{s-t} (d_j-d_0)\prod_{s-t<j} (d_j-d_{s-\ell+1})
+
\prod_{j=s-\ell+1}^{s-t} (d_j-d_0)\prod_{s-t+1<j} (d_j-d_{s-\ell+1})
$$
$$
=
\bigg(\prod_{j=s-\ell+2}^{s-t} (d_j-d_0)\bigg)
\bigg((d_{s-t+1}-d_{s-\ell+1})+(d_{s-\ell+1}-d_0)\bigg)
\bigg(\prod_{s-(t-1)<j} (d_j-d_{s-\ell+1})\bigg)
$$
$$
=\bigg(\prod_{j=s-\ell+2}^{s-(t-1)} (d_j-d_0)\bigg)
\bigg(\prod_{s-(t-1)<j} (d_j-d_{s-\ell+1})\bigg).
$$
\end{proof}

\bigskip

\vbox{\noindent Author Addresses:\par
\smallskip
\noindent{David Eisenbud}\par
\noindent{Department of Mathematics, University of California, Berkeley,
Berkeley CA 94720}\par
\noindent{eisenbud@math.berkeley.edu}\par
\smallskip
\noindent{Frank-Olaf Schreyer}\par
\noindent{Mathematik und Informatik, Universit\"at des Saarlandes, Campus E2 4, 
D-66123 Saarbr\"ucken, Germany}\par
\noindent{schreyer@math.uni-sb.de}\par
}

\end{document}